\documentclass[11pt]{amsart}
\usepackage{graphicx}
\usepackage{lipsum}
\usepackage{multirow}
\usepackage{amsfonts}
\usepackage{subfigure}
\usepackage{graphicx}
\usepackage{epstopdf}
\usepackage{multirow}
\usepackage{color}
\usepackage[ruled,linesnumbered]{algorithm2e}
\usepackage{amsmath}
\newtheorem{theorem}{Theorem}

\newtheorem{remark}{Remark}

\begin{document}
	
\title[An Adam-CBO global optimization method]{A consensus-based global optimization method with adaptive momentum estimation}

\author{Jingrun Chen}
\address{School of Mathematical Sciences and Mathematical Center for Interdisciplinary Research, Soochow University, Suzhou, 215006, China}
\email{jingrunchen@suda.edu.cn}

\author{Shi Jin}
\address{School of Mathematical Sciences, Institute of Natural Sciences, and MOE-LSC, Shanghai Jiao Tong University, Shanghai, 200240, China}
\email{shijin-m@sjtu.edu.cn}

\author{Liyao Lyu}
\address{Department of Computational Mathematics, Science, and Engineering, Michigan State University, East Lansing, MI, 48824, USA}
\email{lyuliyao@msu.edu}

\begin{abstract}
Objective functions in large-scale machine-learning and artificial intelligence applications often live in high dimensions with strong non-convexity and massive local minima. First-order methods, such as the stochastic gradient method and Adam \cite{kingma2014adam}, are often used to find global minima.
Recently, the consensus-based optimization (CBO) method has been introduced as one of the gradient-free optimization methods and its convergence is proven with dimension-dependent parameters, which may suffer from the curse of dimensionality. By replacing the isotropic geometric Brownian motion with the component-wise one, the latest improvement of the CBO method \cite{carrillo2019consensus} is guaranteed to converge to the global minimizer with dimension-independent parameters \cite{ha2019convergence}, although the initial data need to be well-chosen. In this paper, based on the CBO method and Adam, we propose a consensus-based global optimization method with adaptive momentum estimation (Adam-CBO). Advantages of the Adam-CBO method include: (1) capable of finding global minima of non-convex objective functions with high success rates and low costs; (2) can handle non-differentiable activation functions and thus approximate low-regularity functions with better accuracy. The former is verified by approximating the $1000$ dimensional Rastrigin function with $100\%$ success rate at a cost only growing linearly with respect to the dimensionality. The latter is confirmed by solving a machine learning task for partial differential equations with low-regularity solutions where the Adam-CBO method provides better results than the state-of-the-art method Adam. A linear stability analysis is provided to understand the asymptotic behavior of the Adam-CBO method.
\end{abstract}

\date{\today}
\maketitle

\section{Introduction}

The goal of this work is developing consensus-based global optimization methods to solve high dimensional unconstrained optimization problems
\begin{align*}
	x^* = \arg\min_{x\in \mathbb{R}^d} f(x),
\end{align*}
where the target function (loss function) $f(x)$ defined in $\mathbb{R}^d$ achieves a unique global minimizer. 

A high-dimensional nonlinear, non-convex optimization is an essential part of machine learning problems, with the target function defined in general as
\begin{equation*}
	f(x) = \frac{1}{n}\sum_{i=1}^n \|\mathcal{N}_x(\hat{x})-\hat{y}\|,
\end{equation*}
where $x$ is the parameter vector and $\mathcal{N}_x$ represents a neural network representation\footnote{Parameters in a neural network are commonly denoted by $\theta$ instead in Section \ref{subsec: dnn}.}. $(\hat{x}_i,\hat{y}_i)_{i=1}^n$ is a set of labeled data, and $\|\cdot\|$ is the $L^2$ distance between a predicted data point and the corresponding labeled data point.

The gradient descent method, most frequently used method in optimization, often updates the parameters by the iteration scheme 
\begin{equation*}
	x^{t+1} = x^t - \alpha \nabla f(x^t)
\end{equation*}	
with $\alpha$ being the learning rate. However, for a big labeled data set, i.e., $n$ is tremendously big, computing $f$ in each iteration is time consuming, and the iterations often get stuck at local minima. The stochastic gradient descent (SGD) method \cite{bottou2010large,bottou2012stochastic} instead computes $f$ on a randomly selected subset of the labeled data set, by choosing $m$ points randomly from the labeled data set with $m\ll n$ (The subset needs to be updated at each iteration).  The SGD method with momentum term \cite{qian1999Jan} damps oscillations in the SGD method by introducing exponentially weighted moving average as the momentum
\begin{align*}
	&x^{t+1} = x^t - m^t,\\
	&m^t = -\gamma m^{t-1} + \alpha \nabla f(x^t).
\end{align*}
The momentum term increases for dimensions whose gradients point toward the same direction and decreases for dimensions whose gradients change directions. Adding the momentum leads to a faster convergence than the SGD method and shows higher possibility to jump out of local minima. However, if the momentum is added too much, the global minimizer will be most likely missed. The iterator typically rolls past the global minimizer, and then rolls backwards but misses it again. Thus, adding too much momentum often generates a sequence that swings back and forward between local minima. Later, the adaptive momentum method (Adam) \cite{kingma2014adam} also adds the estimation of the second order momentum
\begin{align*}
	&x^{t+1} = x^t - \gamma \frac{\hat{m}^t}{\sqrt{\hat{v}^t}+\epsilon},\\
	&m^t = \beta_1 m^{t-1} + (1-\beta_1) \nabla f(x^t), \quad \hat{m}_t = \frac{m_t}{1-\beta_1^t},\\
	&v^t = \beta_2 v^{t-1} + (1-\beta_2) \nabla^2 f(x^t), \quad \hat{v}_t = \frac{v_t}{1-\beta_1^t},
\end{align*}
where $0<\beta_1,\beta_2<1$. The second order momentum here provides an adaptive adjustment of the learning rate, which has been used in AdaGrad \cite{duchi2011adaptive}, AdaDelta \cite{zeiler2012adadelta}, and RMSprop. By combining the advantages of AdaGrad for dealing with sparse gradients and RMSProp for dealing with non-stationary objectives, the Adam method has been widely used.

However, in many cases the objective function is not differentiable and the training of deep neural networks has the issue of gradient explosion or vanishing \cite{bengio1994learning}. In general, gradient-based methods do not offer a guarantee of global convergence in high dimensional and non-convex problems. Long before machine learning becomes popular, no-convex and nonlinear optimization problems have been considered in some evolutionary computation methods, including the Nelder-Mead method \cite{lagarias1998convergence,nelder1965a}, the genetic algorithm \cite{whitley1994genetic,harik1999compact}, the simulated annealing method \cite{van1987simulated,kirkpatrick1983optimization}, and the particle swarm optimization \cite{poli2007particle,kennedy1995particle}. Despite the tremendous empirical success of these techniques, it is often difficult to provide guarantees of robust convergence to the global minimizer. 
 
The focus of the current work is the CBO method, where a particle system consisting of N particles, labeled as $X^i_t, i = 1,\cdots N$, is considered. During the dynamic evolution, the particle system tends to their weighted average, and meanwhile undergoes some fluctuation due to the random noise, such as the isotropic geometric Brownian motion \cite{poli2007particle,carrillo2018analytical}. Ideally, these particles are expected to gather at the global minimizer of the objective function associated to the system. Mathematically, such a convergence was proved in \cite{carrillo2018analytical} with exponential rate in time under dimension-dependent conditions, i.e., the learning rate depends on the dimension. Therefore, the CBO method may suffer from the curse of dimensionality. To overcome this issue, in \cite{carrillo2019consensus}, Carrillo, Jin, Li, and Zhu proposed to replace the isotropic geometric Brownian motion with the component-wise one. Such a modification leads to the convergence to the global minimizer with dimension-independent parameters, as proved in \cite{ha2019convergence} for well-chosen initial data. From the perspective of efficiency, the idea of random mini-batch is used for quantities involving the summation of individual particle contribution \cite{carrillo2019consensus}, which reduces the computational complexity from $\mathcal{O}(N)$ to $\mathcal{O}(\frac{N}{M})$ with $M$ being the number of particles in each batch. For extremely high dimensional problems, these method require very well-chosen initial data which may be difficult for practical problems.
 
In this work, we improve the CBO method \cite{carrillo2019consensus} by adding first and second order momentum terms to damp the oscillation and accelerate the convergence. In general, we emphasize that the Adam-CBO method has the ability to handle non-differentiable object functions, and has the improved possibility to find the global minimizer of high-dimensional and non-convex functions at a cost only growing linearly with respect to the dimensionality. 
This will be demonstrated by various numerical experiments.

The article is organized by the following structure. In Section \ref{sec:AM-CBO}, we propose the Adam-CBO method together with a brief introduction of the CBO method for completeness. In Section \ref{sec:r-function}, using the example of Rastrigin function, we find that the Adam-CBO method performs better than the CBO method with a higher possibility to find the global minimizer with the same cost. In Section \ref{sec:neural network}, using the Adam-CBO method to approximate functions, we find that the Adam-CBO method also has the spectral bias \cite{Rahaman2018}, or the Frequency principle \cite{xu2019frequency}, which is similar to the first-order methods. In addition, the Adam-CBO method is used to solve partial differential equations (PDEs) with low-regularity solutions. By using activation functions that cannot take gradients, the Adam-CBO outperforms Adam in terms of approximation accuracy. The conclusion is drawn in Section \ref{sec:Conclusion}.

\section{A consensus based optimization method with adaptive momentum estimation}
\label{sec:AM-CBO}

In this section, we provide a detailed discussion on the Adam-CBO method and give some theoretical insights on its convergence.
For completeness, we first give a brief introduction to the CBO method.

\subsection{The CBO method}
\label{subsec:CBO}

The CBO method considers a stochastic interacting system of $N$ particles with position $X^i_t = (x^i_1,\cdots , x^i_d)^T \in \mathbb{R}^d$, whose dynamics can be described as a first order system \cite{carrillo2010asymptotic,carrillo2018analytical, carrillo2019consensus}
\begin{align}
\dot X^i_{t} = -\lambda (X^i_t - x^*)  + \sigma (X_t^i - x^* ) \dot W_t^i, \quad 1\leq i \leq N, \label{equ:CBO model2}
\end{align}
where $\lambda$ represents the learning rate, $N$ is the number of particles, and $M$ is the number of particles in each batch. Here 
\[
x^* = \sum_{i=1}^N X_t^i \frac{\omega_f^\alpha (X^i_t)}{\sum_{j=1}^N\omega^\alpha _f (X^j_t)},
\]
where $\omega_f^\alpha$ is a weight function and can be taken as mode $\omega_f^\alpha = \exp(-\alpha f(x))$ for some appropriately chosen $\alpha >0$, and $f(x)$ is a given (possibly non-convex) function to be optimized. $\dot X$ denotes the temporal derivative of $X$.

We discretize the system \eqref{equ:CBO model2} with stepsize $1$, and obtain 
\begin{align}
\label{equ:CBO model1}
     X^i_{t+1} = X^i_t -\lambda (X^i_t - x^*)  + \sigma (X_t^i - x^*) \mathrm{d}W_t^i,
\end{align}
 The component-wise geometric Brownian motion $W_t^i$ is used to replace the noise in the numerical implementation. Details of the algorithm can be found in \textbf{Algorithm \ref{alg:CBO}}. Without loss of generality, we assume $N//M$. Note that $t_N$ represents the maximum number of temporal steps, or the final time due to the stepsize $1$. If necessary, one can choose a stopping criterion, like $\max_i|X^i_t-x^*| < e $ to stop the update ahead of the final time $t = t_N$.
\begin{algorithm}
  \KwIn{$\lambda$, $N$, $M$, $t_N$}
  \tcc{$\lambda$ represents the learning rate, $N$ is the number of particles, $M$ is the number of particles in each batch and $t_N$ is the number of iterations.}
  Initial $X^i_0$, $i = 1,\cdots N$;

  \For{$t = 0$ \textbf{to} $t_N$ }{
    Generate an index set $P_k$ by random permutation of $\{1,2,\cdots,N\}$;
    
    Generate batch sets of particles in the order of $P_k$ as $B^1,\cdots B^{\frac{N}{M}}$ with each batch having $M$ particles;
    
    \For{$j = 1$  \textbf{to} $\frac{N}{M}$}{
        Update $x^* =\sum\limits_{k\in B^j} \frac{X_t^k\mu^k_t}{\sum\limits_{i\in B^j} \mu_t^i }$, where $\mu_t^i = \omega_f^\alpha (X_t^i)$;
        
        Update $X^i_t$ for $j \in B^j$ as follows
        
        $
        X^i_{t+1} = X^i_t - \lambda \gamma _{k,\theta}(X^i_t-x^*) + \sigma_{k,\theta}\sqrt{\gamma_{k,\theta}} \sum\limits_{k = 1}^d \vec{e}_k(X^i_t-x^*) z_i \quad z_i \sim N(0,1)
        $.
        
        \tcc{$e_k$ is the unit vector along the $k$-th dimension.}
    }
    }
    \KwOut{$X_{t_N}^i, \quad i = 1\cdots N$}
    
  \caption{Consensus-based global optimization method.}
  \label{alg:CBO}
\end{algorithm}

\subsection{The Adam-CBO method}
\label{subsec:Adam-CBO}

By introducing an additional momentum $M^i_t$, we rewrite the first order system \eqref{equ:CBO model2} in Section \ref{subsec:CBO} into
\begin{align}
    \dot X^i_t & = -\lambda M^i_t + \sigma ^t \dot W_t^i, \quad i=1,\cdots,N, \label{eqn:2ndsystem1}\\
    M^i_t & = X^i_t- x^*. \label{eqn:2ndsystem2}
\end{align}
Note that the stochastic term in \eqref{eqn:2ndsystem1} is isotropic since it is found that such a modification leads to a better numerical performance in the Adam-CBO method, while the anisotropic stochastic term in the CBO method performs better with theoretical guarantees \cite{ha2019convergence}.
Discretization of \eqref{eqn:2ndsystem1} yields
\begin{equation}
	X_{t+1}^i= X_{t}^i -\lambda M^i_t + \sigma^t dW_t^i. \label{eqn:adam-cbo1}
\end{equation}
By definition \eqref{eqn:2ndsystem2}, we have
\begin{align*}
M^i_{t+1} & = X^i_{t+1} - x^* = (X^i_{t+1} -X^i_t) + M^i_t \\
& = (1-\lambda) M_t^i + \sigma^t dW_t^i.
\end{align*}
To update the momentum $M_t^i$ adaptively, we borrow the idea from the Adam method \cite{kingma2014adam}. In the asymptotic sense, as $t\rightarrow + \infty$, 
$\sigma^t dW_t^i$ can be represented by $\lambda M_{t+1}^i$. Thus the above equation can be rewritten as
\begin{equation}
M^i_{t+1} = \beta_1 M_t^i + (1-\beta_1) (X_{t+1}^i - x^*) \label{eqn:adam-cbo2}
\end{equation}
with $\beta_1=1-\lambda$.

We now show the relationship between $M^i_t$ and the first moment of $X^i_t-x^*$. Using \eqref{eqn:adam-cbo2} recursively, one gets
\begin{equation*}
	\begin{aligned}
		M^i_t & = \beta_1 M^i_{t-1} + (1-\beta_1)(X^i_t - x^*)\\
		& = \beta_1 (\beta_1 M^i_{t-2} + (1-\beta_1) (X_{t-1}^i - x^*)) + (1-\beta_1) (X^i_t - x^*)\\
		& = \cdots \\
		& = (1-\beta_1) \sum_{k=0}^t \beta_1^{t-k} (X_{k}^i - x^*).
	\end{aligned}
\end{equation*}
Assume that $X_{k}^i - x^*$ is stationary, i.e., they have the same distribution for different $k$, then
\begin{equation*}
	\begin{aligned}
		\mathbb{E} [M^i_t] & =(1-\beta_1)   \mathbb{E} [\sum_{k=0}^t \beta_1^{t-k} (X_{k}^i - x^*)] \\
		&= (1-\beta_1) \mathbb{E}[ X_{t}^i - x^* ]  \sum_{k=0}^t \beta_1^{t-k} \\
		&= (1-\beta_1^t)  \mathbb{E}[ X_{t}^i - x^* ].
	\end{aligned}
\end{equation*}
Therefore, $M^i_t$ gives an estimation of the first moment of $ (X_{k}^i - x^*)$ as $t \rightarrow \infty$. To get an unbiased estimation of $ (X_{k}^i - x^*)$ for small $t$ as well, we rescale $M^i_t$ by $(1-\beta_1^t)$ and denote by $\hat{M}^i_t$ in \textbf{Algorithm \ref{alg:CBO-moment}}. This argument provides a connection between \eqref{eqn:adam-cbo1} and \eqref{equ:CBO model1}.

For the second order moment $\mathbb{E}(|X^i_t-x^*|^2)$\footnote{The square here is defined in the element-wise sense.}, we define 
\begin{equation}
	\begin{aligned}
	\label{equ:CBO second order moment}	V^i_t  = \beta_2 V^i_{t-1} + (1-\beta_2) |X^i_t - x^*|^2. 
	\end{aligned}
\end{equation}
Application of the same argument for $\mathbb{E}[X^i_t]$ yields
\begin{equation}
	\mathbb{E}[V^i_t] = (1-\beta_2^t) \mathbb{E}[|X^i_t-x^*|^2],
\end{equation}
and $\hat{V}^i_t = \frac{V^i_t}{1-\beta_2^t}$ is an unbiased estimation of $\mathbb{E}[|X^i_t- x^*|^2]$. Therefore, we modify \eqref{eqn:adam-cbo1} by 
\begin{equation}
	\begin{aligned}
	\label{equ: CBO moment model true}	X^i_{t+1} = X_t^i - \frac{\lambda \hat{M}^i_{t+1}}{\sqrt{\hat{V}^i_{t+1}} + \epsilon} + \sigma^ t  d W_t^i,
	\end{aligned}
\end{equation}
where $\epsilon$ is a small number and typically takes the value $1e-8$ to avoid the vanishing of the denominator.
Combining \eqref{equ: CBO moment model true}, \eqref{eqn:adam-cbo2}, and \eqref{equ:CBO second order moment} gives \textbf{Algorithm \ref{alg:CBO-moment}}.
Although $\beta_1=1-\lambda$ in the above derivation, $\beta_1$ and $\beta_2$ are chosen to be independent of $\lambda$. In practice, we set $\beta_1 = 0.9$ and $\beta_2=0.99$.

The Adam-CBO method differs from the CBO method in the following aspects. First, it adds estimations of first momentum $M^i_t$ $(\hat{M^i_t})$ and second momentum $V^i_t$ $(\hat{V^i_t})$ into the algorithm without increasing much computational costs. Second, the component-wise geometric Brownian motion term $ \sum\limits_{k = 1}^d \vec{e}_k(X^i_t-x^*) z_i$ is replaced by $\sum_{k = 1}^d \vec{e}_k z_i $, which puts stochastic effects in different dimensions on equal footing. In the case that $X^i_t$ can converge to $x^*$ quickly, so the Adam-CBO method shall have the stronger ability to explore the landscape of the loss function. Note that $\sigma^t$ in \textbf{Algorithm \ref{alg:CBO-moment}} is a decreasing function of $t$, so the method is expected to converge at the finial time. Typically, $\sigma^t = 0.99^{t/10}$ or $\sigma^t = 0.99^{t/100}$ is used in practice.

\begin{algorithm}
\KwIn{$\lambda$, $N$, $M$, $t_N$, $\beta_1$, $\beta_2$}
  \tcc{$\lambda$ represents the learning rate, and $\beta_1,\beta_2$ are the exponential decay rates for the first and the second order moment estimation, respectively.}
  Initialize $X^i_0$, $i = 1,\cdots N$ by the uniform distribution;
  
  Initial $M^i_0,V^i_0 =0$;
  \tcc{Initialize first order and second order moments.}
  
  \For{$t = 0$ \textbf{to} $t_N$ }{
    Generate a random permutation of index $\{1,2,\cdots,N\}$ to form set  $P_k$;
    
    Generate batch set of particles in order of $P_k$ as $B^1,\cdots B^{\frac{N}{M}}$ with each batch having $M$ particles;
    
    \For{$j = 0$  \textbf{to} $\frac{N}{M}$}{
        Update $x^* =\sum\limits_{k\in B^j} \frac{X_t^k\mu^k_t}{\sum\limits_{i\in B^j} \mu_t^i }$, where $\mu_t^i = \omega_f^\alpha (X_t^i)$;
        
        Update $X^i_t$ for $j \in B^j$ as follows
        
        $
        M^i_{t+1} = \beta_1 M_{t}^i +(1-\beta_1) (X^i_t-x^*) \quad \quad \hat{M}^i_{t+1} = M^i_{t+1}/(1-\beta_1^t)
        $;
        
        $
        V^i_{t+1} = \beta_2 V_{t}^i +(1-\beta_2) (X^i_t-x^*)^2 \quad \quad  \hat{V}^i_{t+1} = V^i_{t+1}/(1-\beta_2^t)
        $;
        
        $
        X^i_{t+1} = X^i_t - \lambda \hat{M^i_t}/(\sqrt{\hat{V^i_t}}+\epsilon) + \sigma^t \sum_{k = 1}^d \vec{e}_k z_i  \quad z_i\; \text{ is a random variable}
        $.
    }
    }
    \KwOut{$X_{t_N}^i, \quad i = 1\cdots N$}
  \caption{Consensus-based global optimization method with adaptive momentum estimation.}
  \label{alg:CBO-moment}
\end{algorithm}

Note that the Adam-CBO method is designed to be adaptive by choosing the learning rate (step size) $\lambda$ automatically rather than empirically. It adapts the learning rate to the parameters, and performs smaller updates (low learning rates) for parameters associated with frequently occurring features, and larger updates (high learning rates) for parameters associated with infrequent features. 

\subsection{A linear stability analysis of the Adam-CBO method}
To understand the algorithmic performance, we consider the linearized problems of both methods at the continuous level and prove their convergences. Note that this does not prove the convergence of the Adam-CBO method, but provides an intuitive understanding of it. We first rewrite \textbf{Algorithm \ref{alg:CBO-moment}} into a continuous form and ignore the stochastic term
	\begin{align}
	\label{equ:C1}	\dot m = (\beta_1 -1) m + (1-\beta_1) (x-\bar{x}),\\
	\label{equ:C2}	\dot v = (\beta_2 -1) v + (1-\beta_2) (x-\bar{x})^2,\\
	\label{equ:C3}	\hat{m} = \frac{m}{1-\beta_1^t} \quad  \hat{v} = \frac{v}{1-\beta_2^t}, \\
	\label{equ:C4}	\dot x =  - \lambda \frac{\hat{m}}{\sqrt{\hat{v}}+\epsilon},
	\end{align}
	where $\bar{x}$ is the optimal solution (constant). We shall prove $x \rightarrow \bar{x}$ with a convergence rate independent of $\lambda$ when $x$ is close to $\bar{x}$ by the linear stability analysis. Denote $\tilde{x} = x - \bar{x}$. Linearizing the system \eqref{equ:C1}-\eqref{equ:C4} around $m = 0, x = \bar{x}, v = 0 $, we have
	\begin{align}	
	&\dot m = -(1-\beta_1) m + (1-\beta_1) \tilde{x},\\
	&\dot v = -(1-\beta_2)v, \\
	&\dot {\tilde{x}} = -\frac{\lambda }{(1-\beta_1^t)\epsilon}m \rightarrow -\frac{\lambda}{\epsilon} m = -\mu m  \quad (t \rightarrow \infty)
	\end{align}
	with $\mu=\lambda/\epsilon$, and in a vector form,
	\begin{equation}
	\begin{aligned}
	\partial_t	\left(\begin{matrix}
			m \\ v \\ \tilde x
		\end{matrix}\right) = 
	\left(\begin{matrix}
		-(1-\beta_1) & 0 & 1-\beta_1\\ 
		0 &  -(1-\beta_2) & 0  \\
		 -\mu & 0 & 0  
	\end{matrix}\right) 
	\left(\begin{matrix}
	m \\ v \\ \tilde x
	\end{matrix}\right). 
	\end{aligned}
	\end{equation}

\begin{theorem}
	\textbf{Algorithm \ref{alg:CBO-moment}} generates a sequence that converges to the optimal solution with rates independent of the learning rate $\lambda$. 
	\begin{proof}
	Eigenvalues of the matrix on the right-hand side are $\beta_2-1$ and $\frac{1}{2}( \beta_1 - 1 \pm i  \sqrt{1-\beta_1}\sqrt{\beta_1-1+4\mu})$ (typically $1-\beta_1\ll 4\mu$), respectively. Thus, $m,v,\tilde{x}$ decay to $0$  exponentially with rate $\beta_2 -1 $ when $\beta_1>2\beta_2+1$ and with rate $\frac{1}{2}(\beta_1-1)$ when $\beta_1 < 2 \beta_2 + 1$ in an oscillatory way. 
	\end{proof}
\end{theorem} 

The CBO method without random noise can be written into a continuous form
\begin{equation}
\dot x = - \lambda (x- \bar{x}).
\end{equation}
The ODE can be solved analytically with a decay rate $e^{-\lambda t }$ towards the stationary point.
Therefore, the decay rate of the CBO method depends exponentially on the learning rate $\lambda$.

\begin{remark}
Although the above analysis indicates that the decay rate of the Adam-CBO method is independent of $\lambda$, $\lambda$ does control the oscillatory behavior during the iteration. Therefore, we argue that during the initial training stage, a large $\lambda$ is favored to make particles oscillate and escape local minima. During the finial training stage, to make the particles converge to the global minimizer faster, we often set a smaller $\lambda$ to control the oscillations.  
\end{remark}

\section{The Rastrigin function}
\label{sec:r-function}

\begin{figure}[ht]
	\centering
	\includegraphics[width=0.6\linewidth]{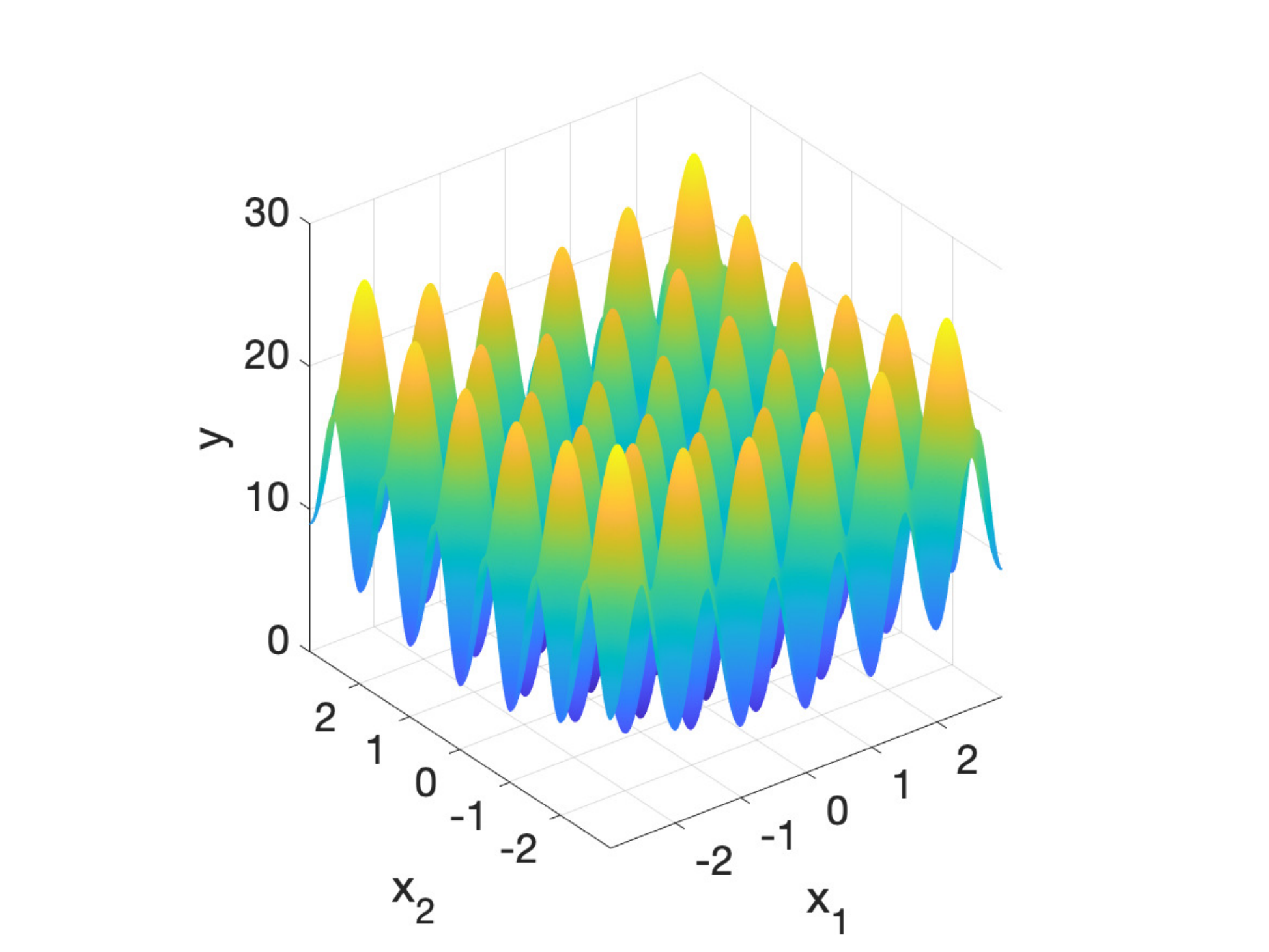}
	\caption{The landscape of the rastrigin function in two dimension with $x\in [-3,3]^2$ and $B = C=0$.}
	\label{fig:Rastrigin}
\end{figure}

In this section, we demonstrate the advantage of the Adam-CBO method by finding the global minimizer of the Rastrigin function
\begin{equation}\label{eqn:rastrigin}
	f(x) = \frac{1}{d} \sum_{i=1}^d \left[(x_i-B)^2 -10\cos(2\pi (x_i-B)) + 10\right] + C
\end{equation}
with $B = \arg \min f(x)$ and $C= \min f(x)$. Figure \ref{fig:Rastrigin} is a visualization of \eqref{eqn:rastrigin} when $d=2$ and $B=C=0$. 

Number of local minima of the rastrigin function in terms of dimension when $B=C=0$ and $x\in(-3,3)$ is listed in Table \ref{tbl:Number of local minimum}.
The number of local minima is $5^d$, which grows exponentially fast in term of the dimensionality. When $d=1000$, the number of minima is $5^{1000}$, approximately $10^{690}$.
\begin{table}[ht]

	\centering\begin{tabular}{|c|c|c|c|c|c|}
	\hline
	d & 1 & 2 & 30 & 100 & 1000 \\
	\hline
	Number of local minima & $5$ & $5^2$ & $5^{30}$ & $5^{100}$ & $5^{1000}$\\
	\hline 
	\end{tabular}
	\caption{Number of local minima for the Rastrigin function in terms of dimension.}
	\label{tbl:Number of local minimum}
\end{table}

Results of CBO and Adam-CBO methods with several random processes, including uniform, Gaussian, and Levy processes are recorded in Table \ref{tbl:process sr}.
In all numerical examples, $\beta_1=0.9$ and $\beta_2=0.99$ in the Adam-CBO method. Here $B$ in \eqref{eqn:rastrigin} is set to be a value between $[-3,3]$ and $X^i_0,\; i=1,\cdots,N$ are initialized between $[-3,3]$. For each case, we run the algorithm $100$ times and check the success rate. It is found that the CBO method fails to find the global minimizer when the dimension is over 30, but the Adam-CBO method still has high success rates even when the dimension reaches $1000$. Moreover, it is found that the Poisson process almost always has higher success rates than uniform and Levy processes.

\begin{table}
	\centering
	\begin{tabular}{|c|c|c|c|c|c|}
		\hline
		\multirow{2}*{$d$}&
		\multirow{2}*{$N$}&
		\multirow{2}*{$M$}&\multicolumn{3}{c|}{CBO}\\
		\cline{4-6}
		~& ~ & ~ &  $\mathcal{N}(0,1)$&$\mathcal{U}(-1,1)$ &Wiener process  \\
		\hline 
		2  & 50 & 40 & 100\% & 100\% & 99\% \\
		10 & 50 & 40 & 100\% & 100\% & 2\%  \\
		20 & 50 & 40 & 98\%  & 22\%  & 0\%  \\
		20 & 50 & 20 & 66\%  & 2\%   & 0\%  \\
		30 & 50 & 40 & 26\%  & 0\%   & 0\%  \\
		30 & 500&5   &  0\%  & 0\%   & 0\%  \\
		\hline
		\multirow{2}*{$d$}&
		\multirow{2}*{$N$}&
		\multirow{2}*{$M$}&\multicolumn{3}{c|}{Adam-CBO}\\
		\cline{4-6}
		~& ~ & ~ &  $\mathcal{N}(0,1)$& $\mathcal{U}(-1,1)$ &Wiener process  \\
		\hline
		30 & 500  & 5  & 99\% & 100\% & 0\% \\
		100& 5000 & 5  & 100\% & 100\% & 0\%\\
		1000& 8000 & 50 & 92\% & 20\% &0\%\\
		\hline
	\end{tabular}
	\caption{Comparison of CBO and Adam-CBO methods with different random processes. Setup of parameters are: $\lambda = 1, \gamma = 0.01, \sigma = 5.1$ in the CBO method with $\mathcal{N}(0,1)$; $\lambda = 0.01, \gamma = 0.1, \sigma = 3$ in the CBO method with $\mathcal{U}(-1,1)$; $\lambda = 0.5, \gamma = 0.1, \sigma = 0.1$ in the CBO method with Wiener process. In the Adam-CBO method, we set $\lambda = 0.1$ and $\sigma^t = 0.99^{\frac{t}{20}}$ in all cases. For each random process, hyper-parameters have been optimized in order to get the best success rate.}
	\label{tbl:process sr}
\end{table}

Next, we compare the dependence of success rates on the batch number of particles in Table \ref{tbl:p_batch sr}. It is observed that the Adam-CBO method usually has higher success rates when the particle batch size $M$ becomes smaller and has higher success rates as the total number of particles grows.
\begin{table}
	\centering
	\begin{tabular}{|c|c|c|c|c|c|c|c|c|}
		\hline
		\multirow{2}*{$d$}&
		\multirow{2}*{$N$}&
		\multirow{2}*{$M$}&\multicolumn{2}{c|}{Adam-CBO}&
		\multirow{2}*{$N$}&
		\multirow{2}*{$M$}&
		\multicolumn{2}{c|}{Adam-CBO}\\
		\cline{4-5}
		\cline{8-9}
		~& ~ & ~ &  $\mathcal{N}(0,1)$& $\mathcal{U}(-1,1)$& ~ & ~ &  $\mathcal{N}(0,1)$& $\mathcal{U}(-1,1)$ \\
		\hline
		100& 1000 & 5   & 87\% & 39\% &5000 & 5 & 100\% & 84\%  \\
		100& 1000 & 10  & 94\% & 60\% &5000 & 10 & 100\% & 100\%   \\
		100& 1000 & 20  & 87\% & 49\% &5000 & 20 & 100\% & 100\%   \\
		100& 1000 & 25  & 77\% & 53\% &5000 & 25 & 100\% & 100\%  \\
		100& 1000 & 50  & 45\% & 8\%  &5000 & 50 & 100\% & 100\%  \\
		100& 1000 & 100 & 2\%  & 0\%  &5000 & 100& 100\% & 100\% \\ \hline
	\end{tabular}
	\caption{Comparison of success rates for different batch numbers when the dimension is $100$, $\lambda = 0.1$, and $\sigma^t = 0.99^{\frac{t}{20}}$.}
	\label{tbl:p_batch sr}
\end{table}
Table \ref{tbl:p_N sr} records the success rate in terms of the number of particles $N$. As $N$ grows, the success rate increases.
\begin{table}
	\centering
	\begin{tabular}{|c|c|c|c|c|}
		\hline
		\multirow{2}*{$d$}&
		\multirow{2}*{$N$}&
		\multirow{2}*{$M$}&\multicolumn{2}{c|}{Adam-CBO}\\
		\cline{4-5}
		~& ~ & ~ &  $\mathcal{N}(0,1)$& $\mathcal{U}(-1,1)$\\
		\hline
		1000& 8000 &  50  & 92\% & 20\% \\
		1000& 10000 & 50  & 100\% & 28\% \\
		1000& 12000 & 50  & 100\% & 28\% \\
		1000& 14000 & 50  & 100\% & 32\% \\
		1000& 16000 & 50  & 100\% & 32\%  \\
		\hline
	\end{tabular}
	\caption{Comparison of success rates for different numbers of particles when the dimension is $1000$, $\lambda = 0.1$, and $\sigma^t= 0.99^{\frac{t}{20}}$.}
	\label{tbl:p_N sr}
\end{table}
One may doubt that the Adam-CBO method shall be sensitive to the initialization. To check this point, instead of choosing initial data randomly, we set initial $X^i_t$ to be $0$, i.e., all particles are initially set to be $0$. Table \ref{tbl:initial sr} shows that the Adam-CBO method still has high success rates.
\begin{table}
	\centering
	\begin{tabular}{|c|c|c|c|c|}
		\hline
		\multirow{2}*{$d$}&
		\multirow{2}*{$N$}&
		\multirow{2}*{$M$}&\multicolumn{2}{c|}{Adam-CBO}\\
		\cline{4-5}
		~& ~ & ~ &  $\mathcal{N}(0,1)$& $\mathcal{U}(-1,1)$\\
		\hline
		30& 500 &  5  & 94\% & 100\% \\
		100& 5000 & 5  & 100\% & 94\% \\
		1000& 10000 & 50  & 100\% & 11\% \\
		\hline
	\end{tabular}
	\caption{Comparison of success rates for different dimensions when $X^i_t$ is initialized by $0$ ($X^i_0 = 0$), $\lambda = 0.1$, and $\sigma^t= 0.99^{\frac{t}{20}}$.}
	\label{tbl:initial sr}
\end{table}

It is worth mentioning that different choices of stochastic terms are used in CBO and Adam-CBO methods. These choices are purely based on numerical experiences. For the Rastrigin function in high dimensions, we observe that the component-wise geometric Brownian motion term $ \sum\limits_{k = 1}^d \vec{e}_k(X^i_t-x^*) z_i$ in \textbf{Algorithm \ref{alg:CBO}} provides better results for the CBO method, while the term $\sum_{k = 1}^d \vec{e}_k z_i $ in \textbf{Algorithm \ref{alg:CBO-moment}} provides better results for the Adam-CBO method. Similar results are observed when applying both methods to neural networks.

\section{Application of the Adam-CBO method on neural networks}
\label{sec:neural network}
In this section we will apply the Adam-CBO method to deep neural networks. For completeness, we briefly introduce deep neural networks (DNNs) and its two applications: approximating functions and solving PDEs. A DNN is constructed by a composition of some basic units which contain activation function $\sigma(x)$ and linear transform $W x + b$. More precisely, we define the simplest network
\begin{equation}
	\mathbb{D}(x;\theta) =  \mathbb N_m( \mathbb N_{m-1}(\cdots \mathbb N_{1} (x))),
\end{equation}
where $\mathbb N_i (x) = \sigma(W^ix + b^i)$. The linear transform $W^i x+b^i$ can transfer a vector $x$ to any dimension, so the output dimension of $\mathbb{N}_i$ can be different. Typically, for $x\in\mathbb{R}^d$, we fix the width by choosing $W^1 \in \mathbb{R}^{n,d}$, $W^{i} \in \mathbb{R}^{n,n}$ for $i = 2,\cdots m-1$, and $W^{m}\in \mathbb{R}^{1,n}$. Therefore, we denote $m$ as the network depth and $n$ as the network width. The parameter set $\theta$ consists of $W^i$ and $b^i$ for $i=1,\cdots,m$, which will be optimized by an optimization method.

For function approximations, we consider to approximate a target function $u(x)$ by a DNN $\mathbb D (x)$ over domain $\Omega$. The objective function is defined as 
\begin{equation}\label{opt:function}
	f(\theta) = \|\mathbb{D}(x;\theta)- u(x)\|^2_{L^2(\Omega)}.
\end{equation}

For a Poisson equation, the target solution $u(x)$ is not given in advance, but it satisfies
\begin{equation}
	\left\{
	\begin{aligned}
	   &-\Delta u = f & x\in \Omega\\
		&u = g & x\in \partial \Omega
	\end{aligned}\right..
\end{equation}
By using the Deep Ritz method \cite{weinan2018deep}, we define the loss function as 
\begin{equation}\label{opt:pde}
	f(\theta) = \int_{\Omega}\frac{1}{2}|\nabla \mathbb{D}(x;\theta) |^2 - f(x)\mathbb{D}(x;\theta) \mathrm{d}x + \eta \int_{\partial\Omega} (\mathbb{D}(x;\theta) - g(x))^2\mathrm{d}x.
\end{equation}
For \eqref{opt:function} and \eqref{opt:pde}, the goal is to find the global minimizer of the following problem
\begin{equation}
	\arg \min_{\theta} f(\theta).
\end{equation}

\subsection{Approximating functions}
\begin{figure}[ht]
	\centering
	\includegraphics[width=0.6\linewidth]{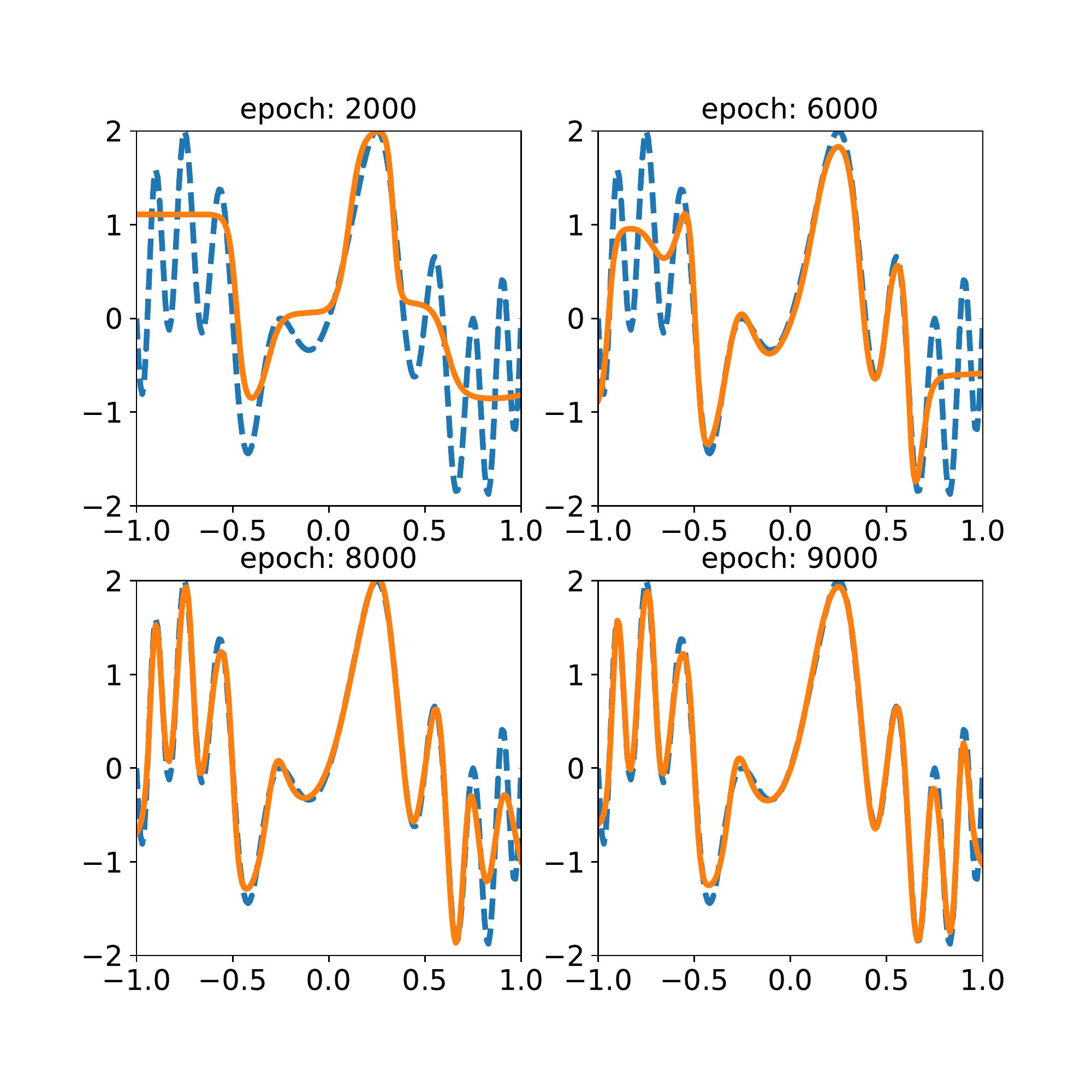}
	\caption{Approximating function \eqref{equ:target 1} using a network with $n = 50$, $m = 3$, and $2701$ parameters in total. The learning rate is $\lambda=0.2$. $N = 500$ particles and $M=5$ particles for each batch are used in the first $50000$ iterations. After that, the random term is ignored and $M = 10$ is used for faster convergence to the optimal solution.}
	\label{fig:fprincipleexm1}
\end{figure}

\begin{figure}[ht]
	\centering
	\includegraphics[width=0.6\linewidth]{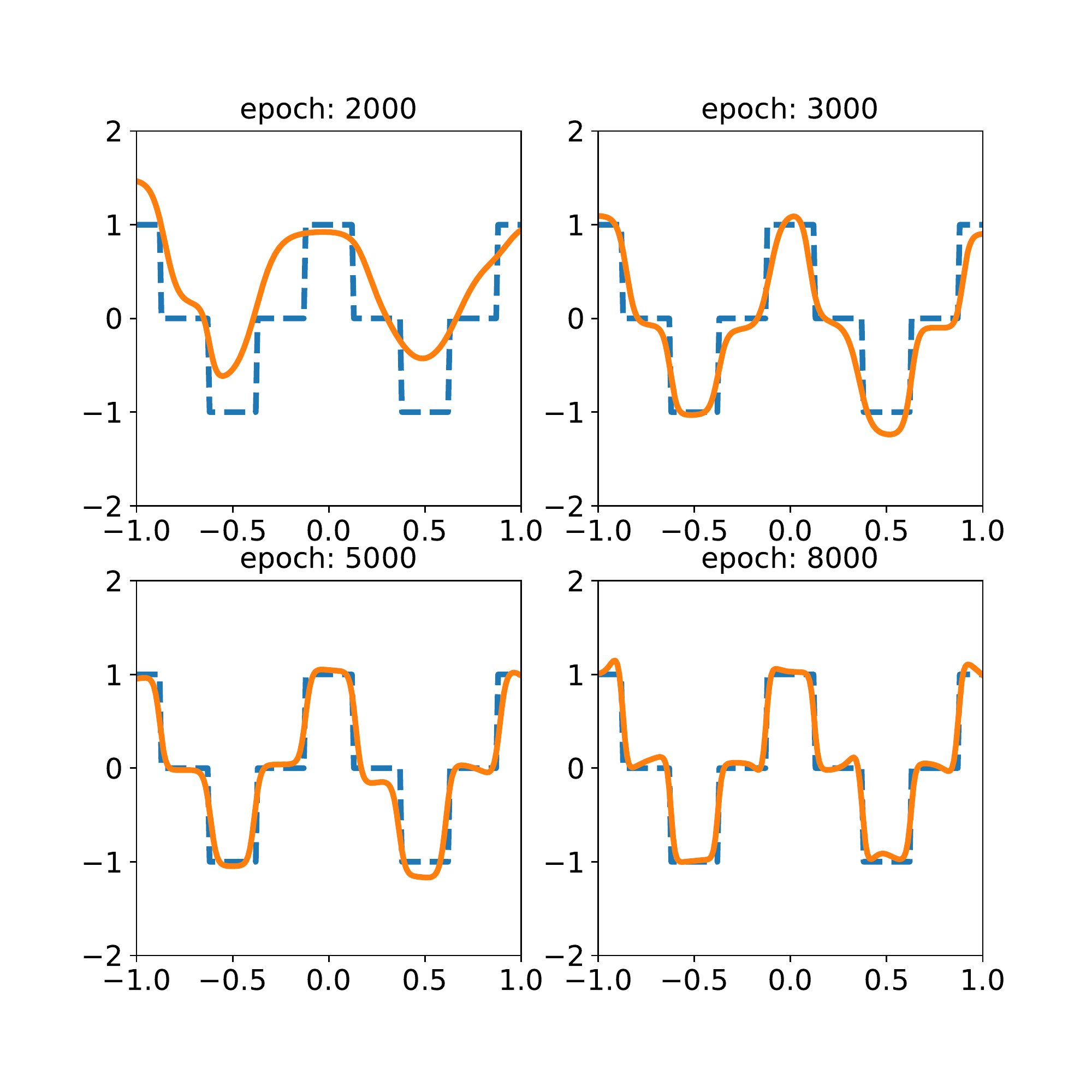}
	\caption{Approximating function \eqref{equ:target 2} using a network with $n = 50$, $m = 3$, and $2701$ parameters in total. The learning rate is $\lambda=0.2$. $N = 500$ particles and $M=5$ particles for each batch are used in the first $50000$ iterations. After that, the random term is ignored and $M = 10$ is used for faster convergence to the optimal solution.}
	\label{fig:fprincipleexm2}
\end{figure}

In this section, we will demonstrate that the Adam-CBO method share the property of spectral bias, or frequency principle as gradient-based methods do \cite{Rahaman2018, xu2019frequency}, i.e., it approximates low-frequency properties of the target function first and high-frequency properties later.
Consider two functions
\begin{align}
	\label{equ:target 1} &u(x) = \sin(2\pi x) + \sin(8 \pi x ^2), \\
	\label{equ:target 2} &u(x) = \left\{\begin{matrix}
		1   & x < -\frac{7}{8}, x> \frac{7}{8}, -\frac{1}{8}<x<\frac{1}{8}\\
		-1  & \frac{3}{8}< x< \frac{5}{8} , -\frac{5}{8}< x<- \frac{3}{8}\\
		0 & \text{otherwise} 
	\end{matrix} \right.,
\end{align}
where the first function is smooth while the second one is not. Here we use the sigmoid function $\sigma(x) = \frac{1}{1+\exp(-x)}$ as the activation function. The training process of the Adam-CBO method is visualized in Figure \ref{fig:fprincipleexm1} for \eqref{equ:target 1} and  in Figure \ref{fig:fprincipleexm2} for \eqref{equ:target 2}. Clearly, the low-frequency information is approximated first and the high-frequency one is captured later. 

\subsection{Deep neural networks}
\label{subsec: dnn}
The commonly used gradient-based method has the issue of gradient vanishing or gradient explosion when the network depth increases.
At the formal level, the Adam-CBO method is independent of the gradient of the loss function with respect to the parameters. Thus it is interesting to check its performance for deeper neural networks. We use DNNs with a fixed width $10$ and different depths to approximate the function
\begin{equation}\label{func3}
	u(x) = \sin(k\pi x ^k).
\end{equation} 
Set $N = 500$ particles, $M = 5$ particles for each batch, and the learning rate $\lambda = 0.2$ in the first $30000$ epochsiterations. Between $30000$ epochs to $80000$ epochs, we set $M = 20$ and ignore the random term to accelerate the convergence. Between $80000$ epochs to $150000$ epochs, we set $M = 100$. After $150000$ epochs, we set the learning rate $\lambda = 1e-2$ to minimize the oscillations. Numerical results are shown in Table \ref{tbl:depth vs accuracy}. For networks with depths $4,7,12,22$, the Adam-CBO method keeps converging to the exact solution. In the implementation, the training process stops after $2\times 10^6$ iterationss. Since the gradient-free method (Adam-CBO) converges slower than the gradient-based method (SGD or Adam), parameters in the neural network fall around the optimal solution but converge to it slowly at the end of the training process. Therefore, it is difficult to obtain the convergence rate of the approximation accuracy in terms of the network depth. However, if SGD or Adam is used with parameters initialized by the uniform distribution, neither method converges well (with final error around $0.3$ in absolute $L^2$ norm) when the network depth is $4$ and $10$, respectively.
\begin{table}
	\centering
	\begin{tabular}{|c|c|c|c|c|}
	\hline
		 depth & Num of parameters & k = 2 &  k = 3 & k = 4  \\
		\hline
		4   & 141 & 6.62 e-03 & 1.32 e-02 & 1.71 e-01\\
		7   & 471 & 4.78 e-03 & 1.42 e-02 & 7.54 e-03\\
		12  & 1021 & 7.44 e-03 & 1.30 e-02 & 5.32 e-02\\
		22  & 2121 & 1.00 e-02 & 1.01 e-02 & 1.21 e-01\\
		\hline
	\end{tabular}
	\caption{Dependence of approximation error measured in absolute $L^2$ norm in terms of network depth for \eqref{func3} when $k=2, 3, 4$.}
	\label{tbl:depth vs accuracy}
\end{table}

\subsection{Solving PDEs with low-regularity solutions}

In this section, we will use the Adam-CBO method to solve PDEs with low-regularity solutions. There has an increasing interest in the development of machine-learning method for solving PDEs; see \cite{weinan2020} for review and references therein. For the purpose of low-regularity solutions, we adopt the Deep Ritz method (DRM) \cite{weinan2018deep}, which is based on the variational formulation associated to the PDE. Consider an elliptic PDE
\begin{equation}\label{eqn:pde}
\left\{
\begin{aligned}
&-\nabla \cdot ( A(x) \nabla u) = - \sum_{i=1}^d\delta(x_i) & x\in \Omega=[-1,1]^d\\
&u(x) = g(x)  & x\in \partial \Omega
\end{aligned}\right.
\end{equation} 
with
\begin{equation}
	A(x)= 
	\left[\begin{matrix}
		(x_1^2)^{\frac{1}{4}} & & \\
		& \ddots & &\\
		& &  (x_d^2)^{\frac{1}{4}}
	\end{matrix}\right].
\end{equation}
The exact solution $u(x)= \sum_{i=1}^d|x_i|^{\frac{1}{2}}$. One can see that the solution is only in $H^{1/2}(\Omega)$ and has singularities when evaluating its derivative at $x_i =0$. The loss function in DRM reads as 
\begin{equation}
	I[u] = \int_{\Omega}\frac{1}{2}(\nabla u)^T  A(x)  \nabla u(x)\mathrm{d}x + \sum_{i=1}^d\int_{-1}^{1}\delta(x_i)u(x)\mathrm{d}x_i + \eta \int_{\partial \Omega} (u(x)-g(x))^2 \mathrm{d}x,
\end{equation}
where $\eta=500$ is the penalty parameter for the boundary condition.

\begin{figure}
	\subfigure[$L^{\infty}$ error]{
		\includegraphics[width=0.45\linewidth]{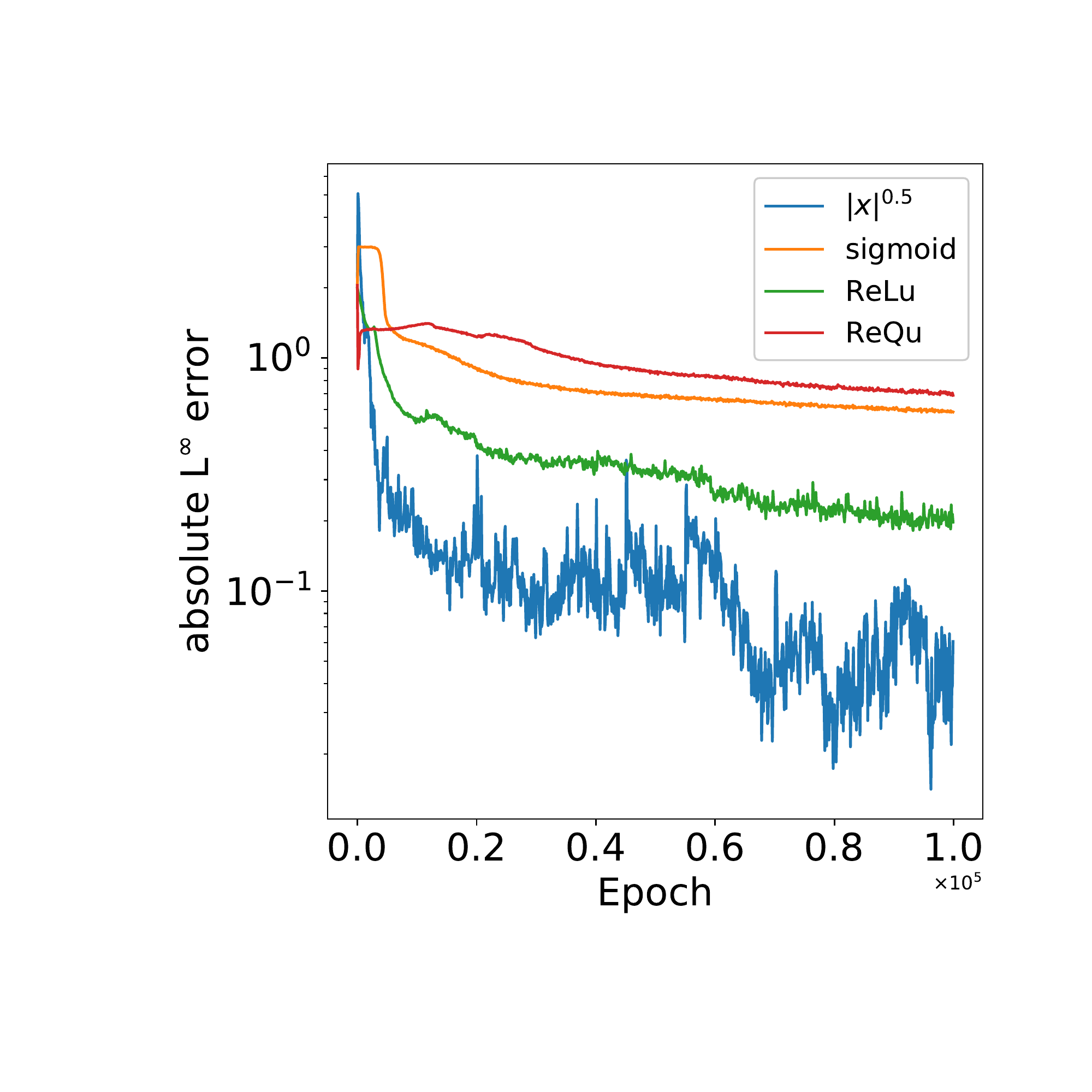}
	}
	\subfigure[$L^2$ error]{
		\includegraphics[width=0.45\linewidth]{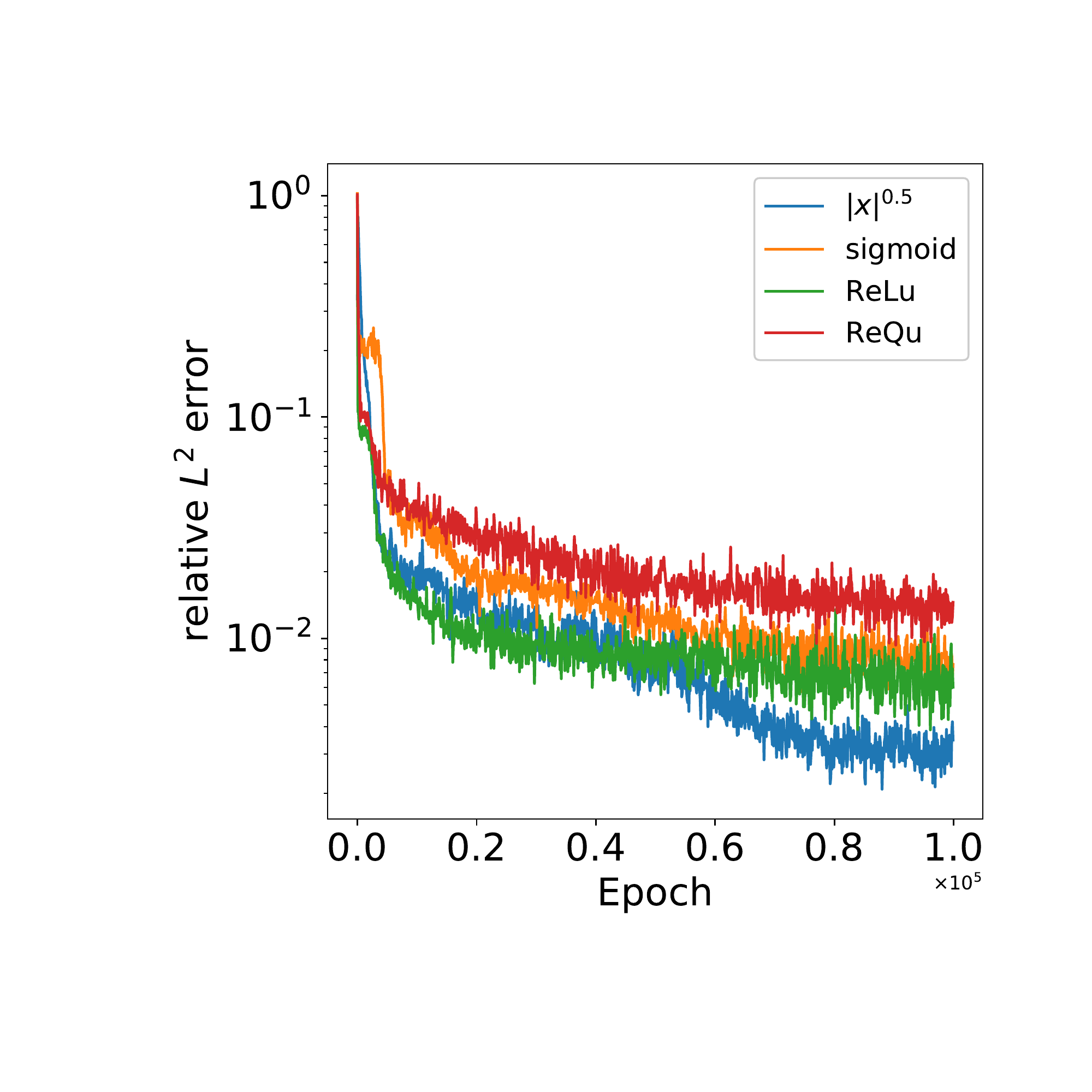}
	}
	\caption{Training process of Adam and Adam-CBO methods for \eqref{eqn:pde} when the dimension is $4$. (a) $L^{\infty}$ error; (b) $L^2$ error.}
	\label{fig:training process singular problem}
\end{figure}

\begin{figure}
	\subfigure[$x_2=x_3=x_4=0$]{
		\includegraphics[width=0.45\linewidth]{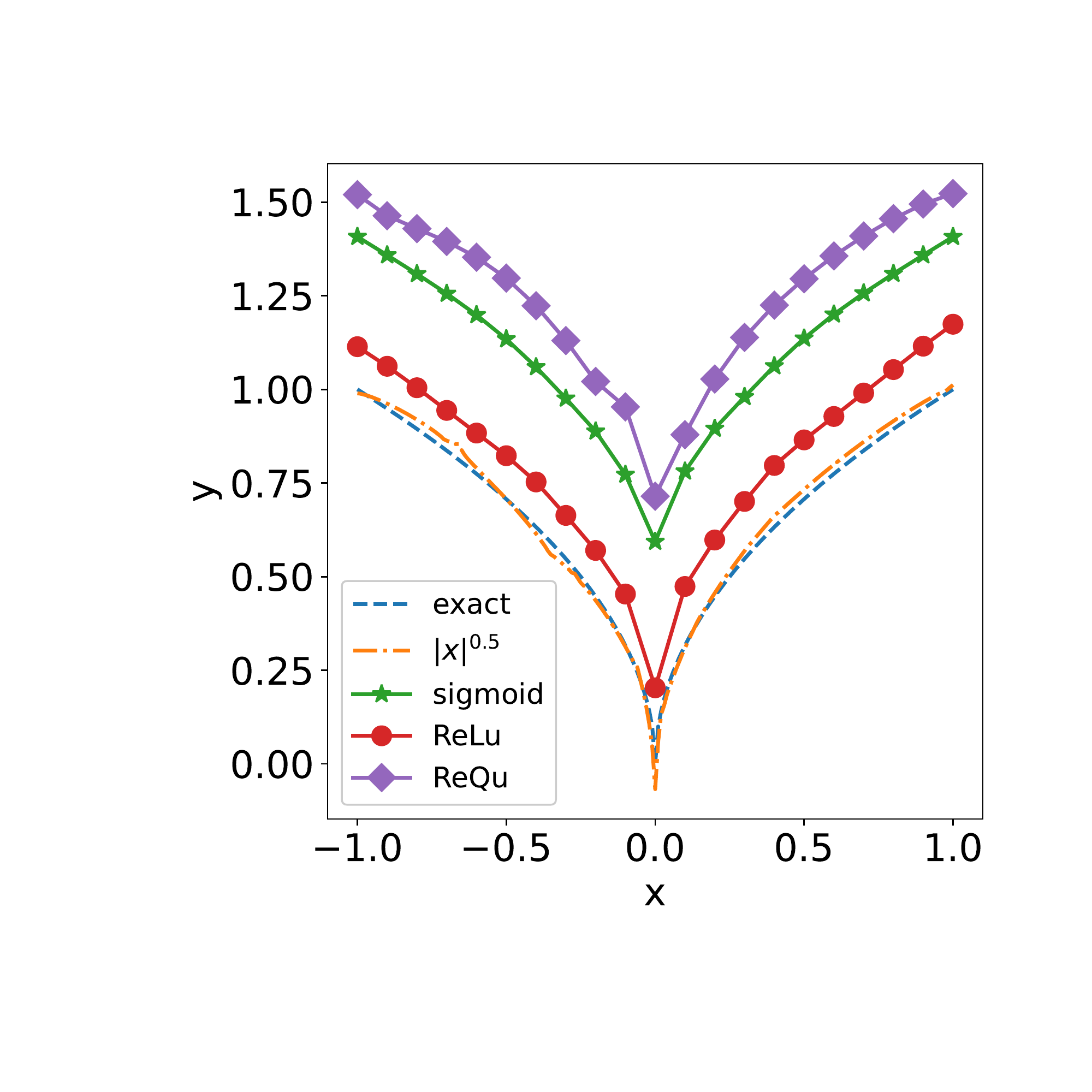}
	}
	\subfigure[$x_1=x_3=x_4=0$]{
		\includegraphics[width=0.45\linewidth]{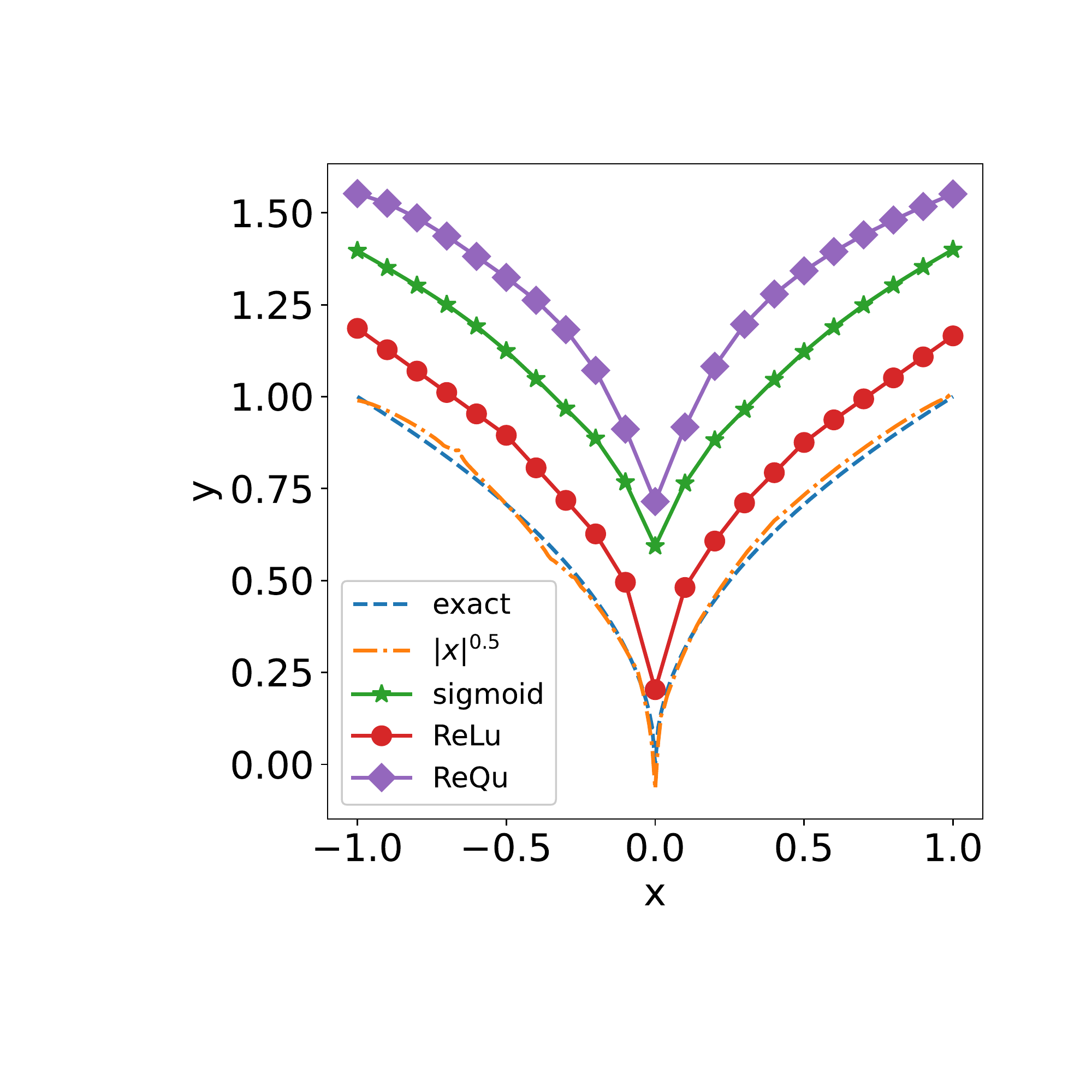}
	}
	\subfigure[$x_1=x_2=x_4=0$]{
		\includegraphics[width=0.45\linewidth]{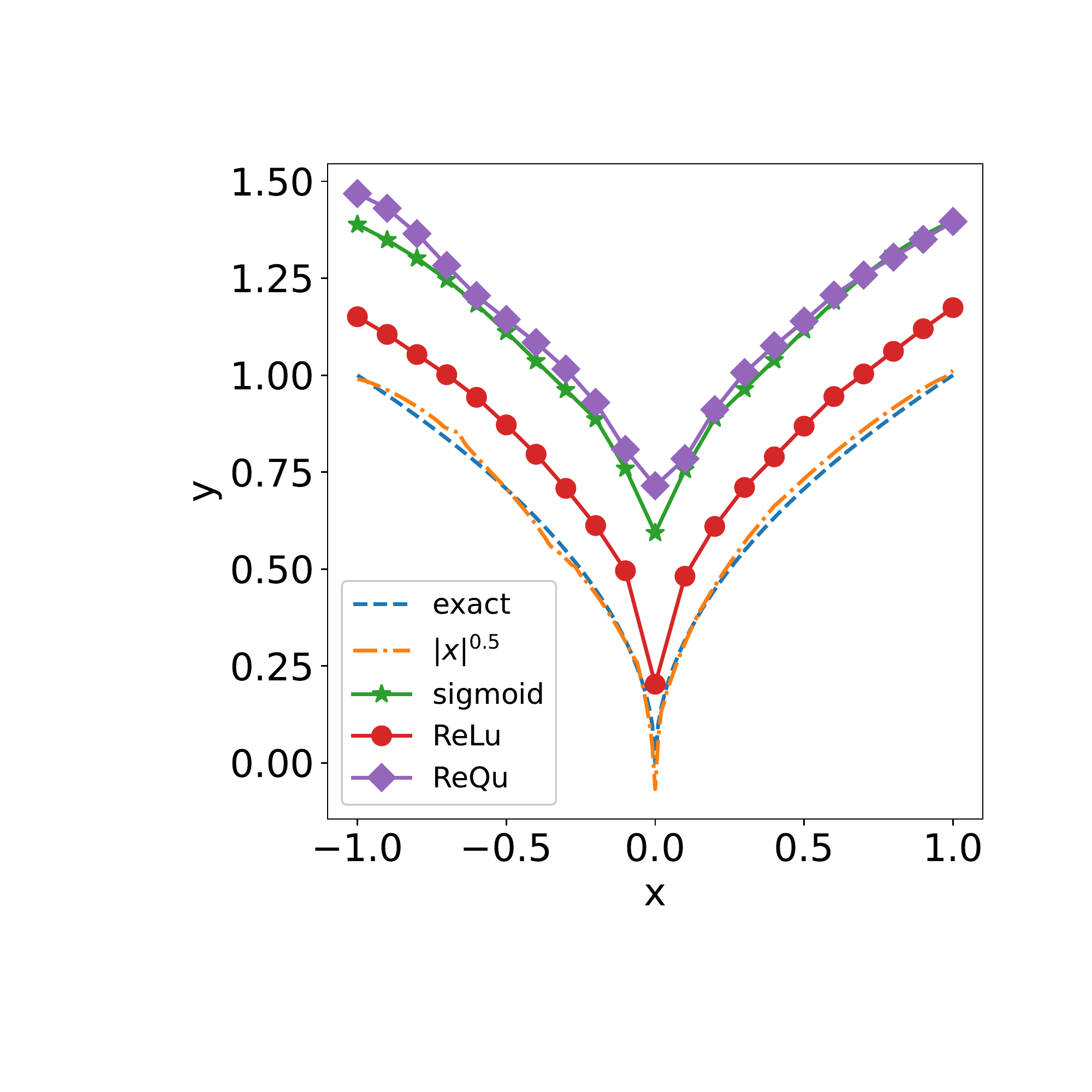}
	}
	\subfigure[$x_1=x_2=x_3=0$]{
		\includegraphics[width=0.45\linewidth]{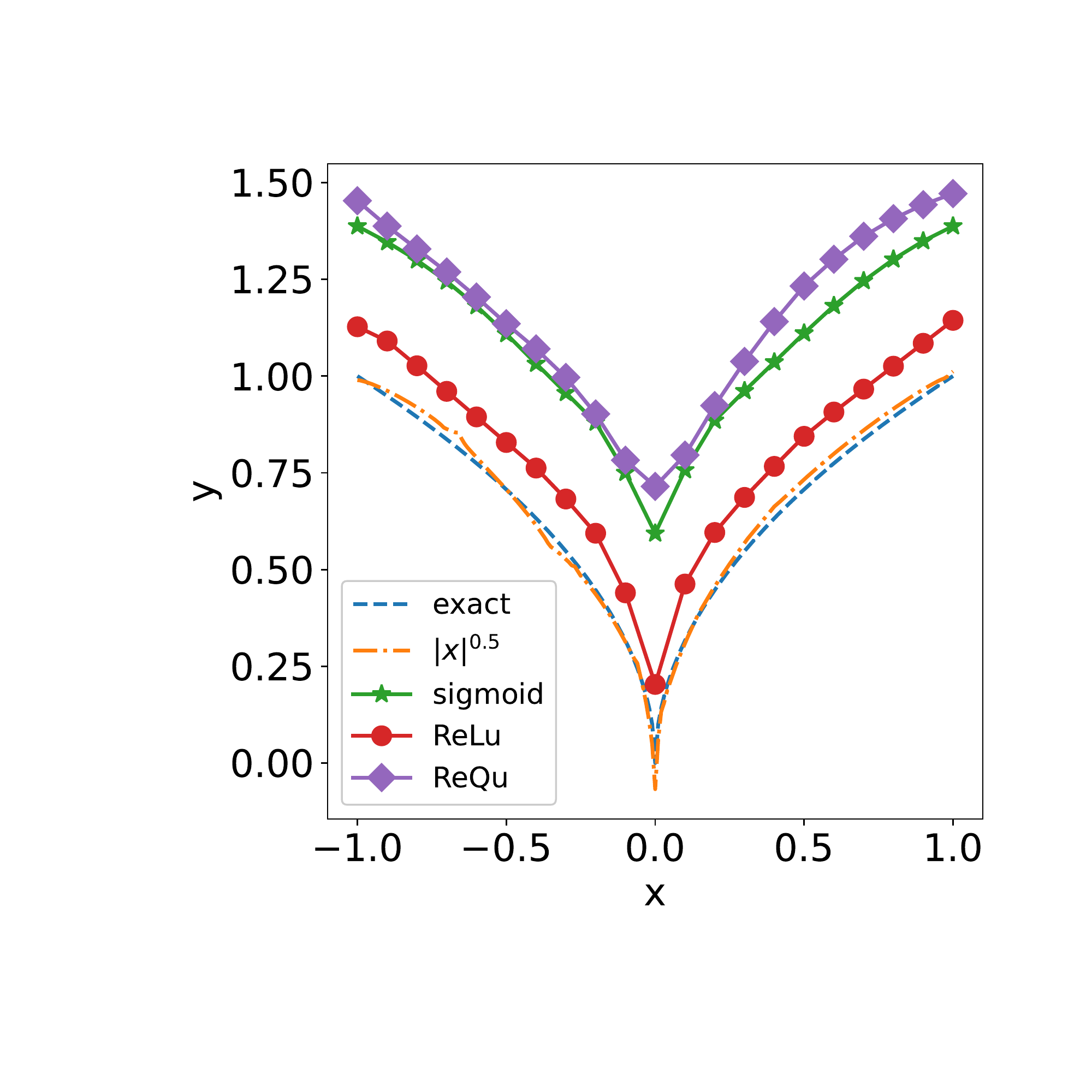}
	}
	\caption{One-dimensional solution profiles at the intersection. (a) $x_2=x_3=x_4=0$; (b) $x_1=x_3=x_4=0$; (c) $x_1=x_2=x_4=0$; (d) $x_1=x_2=x_3=0$.}
\label{fig:landscape of solution of singular problem}
\end{figure}

Activation functions used in Adam include ReLu ($\max\{x,0\}$), ReQu ($\left(\max\{x,0\}\right)^2$), and sigmoid ($\frac{1}{1+\exp(-x)}$). Since the Adam-CBO method is a gradient-free method, to demonstrate its advantage, we use $|x|^{\frac{1}{2}}$ as the activation function. Another reason for choosing this activation function is its low regularity, which leads to superior approximation accuracy in this case. Note that a loss function including this activation function is not differentiable and thus gradient-based methods are not applicable. Numerical results are shown in Table \ref{tbl:error_singular_PDE}. The training process is shown in Figure \ref{fig:training process singular problem}. One-dimensional solution profiles at the intersection where other coordinates are set to be $0$ are visualized in Figure \ref{fig:landscape of solution of singular problem}. It is found that the Adam-CBO method provides better results than Adam with different activation functions. This attributes to the usage of non-differentiable activation functions which better approximate low-regularity PDEs. Moreover, the $|x|^{0.5}$ activation function approximates the low-regularity solution better than any other activation functions near the singularity. 
\begin{table}
	\centering
	\begin{tabular}{|c|c|c|c|c|c|}
		\hline
		d & n & m &  Activation-Optimizer & $L^2$ error & $L^{\infty}$ error\\
		\hline
		\multirow{4}*{2}&\multirow{4}*{20} & \multirow{4}*{2} & ReLu-Adam & 1.23 e-02 & 9.91 e-02\\
		& &  & ReQu-Adam & 2.22 e-02  & 4.21 e-01 \\
		& &  & sigmoid-Adam & 2.19 e-02 & 3.14 e-01\\
		& &  & $|x|^{0.5}$ - Adam-CBO & 3.96 e-03 & 2.09 e-02\\
		\hline
		\multirow{4}*{4}&\multirow{4}*{40} & \multirow{4}*{2} & ReLu-Adam & 6.72 e-03 & 3.70 e-01\\
		& &  & ReQu-Adam & 1.43 e-02 & 1.10 e -00\\
		& &  & sigmoid-Adam & 7.90 e-03 & 7.66 e -02\\
		& &  & $|x|^{0.5}$ -Adam-CBO & 3.13 e-03 & 9.52 e -02\\
		\hline
	\end{tabular}
\caption{Errors measured in $L^2$ and $L^{\infty}$ norms for \eqref{eqn:pde} by Adam and Adam-CBO methods.}\label{tbl:error_singular_PDE}
\end{table}

\section{Conclustion}
\label{sec:Conclusion}

In this work, we propose a consensus-based global optimization method with adaptive momentum estimation based on the consensus-based global optimization method and the adaptive momentum estimation. It shows strong abilities to find global minima for high dimensional problems, including given functions in high dimensions and approximation of low-regularity solutions to PDEs by deep neural networks. The computational complexity is found to  grow linearly with respect to the dimension of the parameter space. Since it is free of gradient, the Adam-CBO method is a suitable choice for problems where derivatives with respect to parameters do not exist. Therefore, it will be of great interest to find the application of the Adam-CBO method for machine learning tasks where non-differentiable activation functions are needed and the dimensionality of parameter space is high.

\medskip
\noindent \textbf{Acknowledgment.} This work of J. Chen was supported by National Key R\&D Program of China (No. 2018YFB0204404) and NSFC grant No. 11971021. The work of S. Jin was supported by NSFC grant No. 11871297.

\bibliographystyle{plain}
\bibliography{references}
\end{document}